\documentclass[runningheads]{llncs}

\usepackage[utf8]{inputenc}
\usepackage{amsmath}
\usepackage{amssymb}
\usepackage{stmaryrd}
\usepackage{multirow}

\newcommand{\R}{\mathbb{R}}
\newcommand{\M}{\mathcal{M}}
\newcommand{\mc}{\mathcal}
\newcommand{\mf}{\mathfrak}
\newcommand{\lto}{\longrightarrow}
\newcommand{\lmto}{\longmapsto}

\newcommand{\ls}{\leqslant}
\newcommand{\gs}{\geqslant}

\newcommand{\tr}{\mathrm{tr}}

\newcommand{\inv}{\mathrm{inv}}

\newcommand{\pow}{\mathrm{pow}}
\newcommand{\diag}{\mathrm{diag}}
\newcommand{\Diag}{\mathrm{Diag}}
\newcommand{\Diff}{\mathrm{Diff}}
\newcommand{\adj}{\mathrm{adj}}
\newcommand{\pol}{\mathrm{pol}}
\newcommand{\Id}{\mathrm{Id}}

\newcommand{\Stab}{\mathrm{Stab}}
\newcommand{\Sym}{\mathrm{Sym}}
\newcommand{\Cov}{\mathrm{Cov}}
\newcommand{\Vect}{\mathrm{span}}
\newcommand{\dotprod}[2]{\langle #1|#2\rangle}

\newcommand{\fun}[4]{
\left\{
\begin{array}{ccc}
#1 & \lto & #2\\ \relax
#3 & \lmto & #4\\
\end{array}
\right.
}

\begin{document}

\title{Is affine-invariance well defined on SPD matrices? A principled continuum of metrics}
\titlerunning{Power-affine and deformed-affine metrics on SPD matrices}

\author{Yann Thanwerdas \and Xavier Pennec}
\authorrunning{Y. Thanwerdas, X. Pennec}

\institute{Universit\'e C\^ote d'Azur, Inria, Epione, France}

\maketitle

\begin{abstract}
Symmetric Positive Definite (SPD) matrices have been widely used in medical data analysis and a number of different Riemannian metrics were proposed to compute with them. However, there are very few methodological principles guiding the choice of one particular metric for a given application. Invariance under the action of the affine transformations was suggested as a principle. Another concept is based on symmetries. However, the affine-invariant metric and the recently proposed polar-affine metric are both invariant and symmetric. Comparing these two cousin metrics leads us to introduce much wider families: power-affine and deformed-affine metrics. Within this continuum, we  investigate other principles to restrict the family size.

\end{abstract}

\section{Introduction}

Symmetric positive definite (SPD) matrices have been used in many different contexts. In diffusion tensor imaging for instance, a diffusion tensor is a 3-dimensional SPD matrix \cite{lenglet_statistics_2006,pennec_riemannian_2006,Fletcher07}; in brain-computer interfaces (BCI) \cite{Barachant13}, in functional MRI \cite{Deligianni11} or in computer vision \cite{Cheng13}, an SPD matrix can represent a covariance matrix of a feature vector, for example a spatial covariance of electrodes or a temporal covariance of signals in BCI. In order to make statistical operations on SPD matrices like interpolations, computing the mean 
or performing a principal component analysis, it has been proposed to consider the set of SPD matrices as a manifold and to provide it with some geometric structures like a \textit{Riemannian metric}, a transitive \textit{group action} or some \textit{symmetries}. These structures can be more or less natural depending on the context of the applications, and they can provide closed-form formulas and consistent algorithms \cite{pennec_riemannian_2006,Dryden09}.

Many Riemannian structures have been introduced over the manifold of SPD matrices \cite{Dryden09}: Euclidean, log-Euclidean, affine-invariant, Cholesky, square root, power-Euclidean, Procrustes... Each of them has different mathematical properties that can fit the data in some problems but can be inappropriate in some other contexts: for example the curvature can be null, positive, negative, constant, not constant, covariantly constant... These properties on the curvature have some important consequences on the way we interpolate two points, on the consistence of algorithms, and more generally on every statistical operation one could want to do with SPD matrices. Therefore, a natural question one can ask is: given the practical context of an application, how should one choose the metric on SPD matrices? Are there some relations between the mathematical properties of the geometric structure and the intrinsic properties of the data?

In this context, the affine-invariant metric \cite{pennec_riemannian_2006,Lenglet06,Fletcher07} was introduced to give an invariant computing framework under affine transformations of the variables. This metric endows the manifold of SPD matrices with a structure of a Riemannian symmetric space. Such spaces have a covariantly constant curvature, thus they share some convenient properties with constant curvature spaces but with less constraints. It was actually shown that there exists not only one but a one-parameter family that is invariant under these affine transformations \cite{Pennec08}. More recently, \cite{Su12,Su14,Zhang18} introduced another Riemannian symmetric structure that does not belong to the previous one-parameter family: the polar-affine metric.

In this work, we unify these two frameworks by showing that the polar-affine metric is a square deformation of the affine-invariant metric (Section 2). We generalize in Section 3.1 this construction to a family of power-affine metrics that comprises the two previous metrics, and in Section 3.2 to the wider family of deformed-affine metrics. Finally, we propose in Section 4 a theoretical approach in the choice of subfamilies of the deformed-affine metrics with relevant properties.

\enlargethispage{8mm}

\section{Affine-invariant versus polar-affine}

The affine-invariant metric \cite{pennec_riemannian_2006,Lenglet06,Fletcher07} and the polar-affine metric \cite{Zhang18} are different but they both provide a Riemannian symmetric structure to the manifold of SPD matrices. Moreover, both claim to be very naturally introduced. The former uses only the action of the real general linear group ${GL}_n$ on covariance matrices. The latter uses the canonical left action of ${GL}_n$ on the left coset space ${GL}_n/O_n$ and the polar decomposition ${GL}_n\simeq{SPD}_n\times O_n$, where $O_n$ is the orthogonal group. Furthermore, the affine-invariant framework is exhaustive in the sense that it provides \textit{all} the metrics invariant under the chosen action \cite{Pennec08} whereas the polar-affine framework only provides \textit{one} invariant metric.

In this work, we show that the two frameworks coincide on the same quotient manifold ${GL}_n/O_n$ but differ because of the choice of the diffeomorphism between this quotient and the manifold of SPD matrices. In particular, we show that there exists a one-parameter family of polar-affine metrics and that any polar-affine metric is a square deformation of an affine-invariant metric. 

In 2.1 and 2.2, we build the affine-invariant metrics $g^1$ and the polar-affine metric $g^2$ in a unified way, using indexes $i\in\{1,2\}$ to differentiate them. First, we give explicitly the action $\eta^i:{GL}_n\times{SPD}_n\lto{SPD}_n$ and the quotient diffeomorphism $\tau^i:{GL}_n/O_n\lto{SPD}_n$; then, we explain the construction of the orthogonal-invariant scalar product $g^i_{I_n}$ that characterizes the metric $g^i$; finally, we give the expression of the metrics $g^1$ and $g^2$. In 2.3, we summarize the results and we focus on the Riemannian symmetric structures of ${SPD}_n$.



\subsection{The one-parameter family of affine-invariant metrics}


\subsubsection{Affine action and quotient diffeomorphism}
In many applications, one would like the analysis of covariance matrices to be invariant under affine transformations $X\lmto AX+B$ of the random vector $X\in\R^n$, where $A\in{GL}_n$ and $B\in\R^n$. Then the covariance matrix $\Sigma=\Cov(X)$, is modified under the transformation $\Sigma\lmto A\Sigma A^\top$. This transformation can be thought as a transitive Lie group action $\eta^1$ of the general linear group on the manifold of SPD matrices:
\begin{equation}
    \eta^1:\fun{{GL}_n\times{SPD}_n}{{SPD}_n}{(A,\Sigma)}{\eta^1_A(\Sigma)=A\Sigma A^\top}.
\end{equation}
This transitive action induces a diffeomorphism between the manifold ${SPD}_n$ and the quotient of the acting group ${GL}_n$ by the stabilizing group $\Stab^1(\Sigma)=\{A\in{GL}_n,\eta^1(A,\Sigma)=\Sigma\}$ at any point $\Sigma$. It reduces to the orthogonal group at $\Sigma=I_n$ so we get the quotient diffeomorphism $\tau^1$:
\begin{equation}
    \tau^1:\fun{{GL}_n/O_n}{{SPD}_n}{[A]=A.O_n}{\eta^1(A,I_n)=AA^\top}.
\end{equation}

\subsubsection{Orthogonal-invariant scalar product}
We want to endow the manifold $\M={SPD}_n$ with a metric $g^1$ invariant under the affine action $\eta^1$, i.e. an affine-invariant metric. As the action is transitive, the metric at any point $\Sigma$ is characterized by the metric at one given point $I_n$. As the metric is affine-invariant, this scalar product $g_{I_n}$ has to be invariant under the stabilizing group of $I_n$. As a consequence, the metric $g^1$ is characterized by a scalar product $g^1_{I_n}$ on the tangent space $T_{I_n}\M$ that is invariant under the action of the orthogonal group.

The tangent space $T_{I_n}\M$ is canonically identified with the vector space $\Sym_n$ of symmetric matrices by the differential of the canonical embedding $\M\hookrightarrow\Sym_n$. Thus we are now looking for all the scalar products on symmetric matrices that are invariant under the orthogonal group. Such scalar products are given by the following formula \cite{Pennec08}, where $\alpha>0$ and $\beta>-\frac{\alpha}{n}$: for all tangent vectors $V_1,W_1\in T_{I_n}\M$, $g^1_{I_n}(V_1,W_1)=\alpha\,\tr(V_1W_1)+\beta\,\tr(V_1)\tr(W_1)$.

\subsubsection{Affine-invariant metrics}
To give the expression of the metric, we need a linear isomorphism between the tangent space $T_\Sigma\M$ at any point $\Sigma$ and the tangent space $T_{I_n}\M$. Since the action $\eta^1_{\Sigma^{-1/2}}$ sends $\Sigma$ to $I_n$, its differential given by $T_\Sigma\eta^1_{\Sigma^{-1/2}}:V\in T_\Sigma\M\lmto V_1=\Sigma^{-1/2}V\Sigma^{-1/2}\in T_{I_n}\M$ is such a linear isomorphism.
Combining this transformation with the expression of the metric at $I_n$ and reordering the terms in the trace, we get the general expression of the affine-invariant metric: for all tangent vectors $V,W\in T_\Sigma\M$,
\begin{equation}
    g^1_\Sigma(V,W)=\alpha\,\tr(\Sigma^{-1}V\Sigma^{-1}W)+\beta\,\tr(\Sigma^{-1}V)\tr(\Sigma^{-1}W).
\end{equation}

As the geometry of the manifold is not much affected by a scalar multiplication of the metric, we often drop the parameter $\alpha$, as if it were equal to 1, and we consider that this is a one-parameter family indexed by $\beta>-\frac{1}{n}$.


\subsection{The polar-affine metric}

\subsubsection{Quotient diffeomorphism and affine action}
In \cite{Zhang18}, instead of defining a metric directly on the manifold of SPD matrices, a metric is defined on the left coset space ${GL}_n/O_n=\{[A]=A.O_n,A\in{GL}_n\}$, on which the general linear group ${GL}_n$ naturally acts by the left action $\eta^0:(A,[A'])\lmto[AA']$. Then this metric is pushed forward on the manifold ${SPD}_n$ into the polar-affine metric $g^2$ thanks to the polar decomposition $\pol:A\in{GL}_n\lmto(\sqrt{AA^\top},{\sqrt{AA^\top}}^{-1}A)\in{SPD}_n\times O_n$ or more precisely by the quotient diffeomorphism $\tau^2$:
\begin{equation}
    \tau^2:\fun{{GL}_n/O_n}{{SPD}_n}{A.O_n}{\sqrt{AA^\top}}.
\end{equation}
This quotient diffeomorphism induces an action of the general linear group ${GL}_n$ on the manifold ${SPD}_n$, under which the polar-affine metric will be invariant:
\begin{equation}
    \eta^2:\fun{{GL}_n\times{SPD}_n}{{SPD}_n}{(A,\Sigma)}{\eta^2_A(\Sigma)=\sqrt{A\Sigma^2A^\top}}.
\end{equation}
It is characterized by $\eta^2(A,\tau^2(A'.O_n))=\tau^2(\eta^0(A,A'.O_n))$ for $A,A'\in{GL}_n$.

\subsubsection{Orthogonal-invariant scalar product}
The polar-affine metric $g^2$ is characterized by the scalar product $g^2_{I_n}$ on the tangent space $T_{I_n}\M$. This scalar product is obtained by pushforward of a scalar product $g^0_{[I_n]}$ on the tangent space $T_{[I_n]}({GL}_n/O_n)$. It is itself induced by the Frobenius scalar product on $\mf{gl}_n=T_{I_n}{GL}_n$, defined by $\dotprod{v}{w}_\mathrm{Frob}=\tr(vw^\top)$, which is orthogonal-invariant. This is summarized on the following diagram.
$$
\begin{array}{ccccc}
{GL}_n & \overset{s}{\lto} & {GL}_n/O_n & \overset{\tau^2}{\lto} & \M={SPD}_n\\ \relax
A & \lmto & A.O_n & \lmto & \sqrt{AA^\top}\\ \relax
\dotprod{\cdot}{\cdot}_\mathrm{Frob} & & g^0_{[I_n]} & & g^2_{I_n}
\end{array}
$$
Finally, we get the scalar product $g^2_{I_n}(V_2,W_2)=\tr(V_2W_2)$ for $V_2,W_2\in T_{I_n}\M$.

\subsubsection{Polar-affine metric}
Since the action $\eta^2_{\Sigma^{-1}}$ sends $\Sigma$ to $I_n$, a linear isomorphism between tangent spaces is given by the differential of the action 
$T_\Sigma\eta^2_{\Sigma^{-1}}:V\in T_\Sigma M\lto V_2=\Sigma^{-1}T_\Sigma\pow_2(V)\Sigma^{-1}\in T_{I_n}\M$. Combined with the above expression of the scalar product at $I_n$, we get the following expression for the polar affine metric: for all tangent vectors $V,W\in T_\Sigma\M$,
\begin{equation}
    g^2_\Sigma(V,W)=\tr(\Sigma^{-2}\,T_\Sigma\pow_2(V)\,\Sigma^{-2}\,T_\Sigma\pow_2(W)).
\end{equation}

\subsection{The underlying Riemannian symmetric manifold}

In the affine-invariant framework, we started from defining the affine action $\eta^1$ (on covariance matrices) and we inferred the quotient diffeomorphism $\tau^1:({GL}_n/O_n,\eta^0)\lto({SPD}_n,\eta^1)$. In the polar-affine framework, we started from defining the quotient diffeomorphism $\tau^2:{GL}_n/O_n\lto{SPD}_n$ (corresponding to the polar decomposition) and we inferred the affine action $\eta^2$. The two actually correspond to the same underlying affine action $\eta^0$ on the quotient ${GL}_n/O_n$. Then there is also a one-parameter family of affine-invariant metrics on the quotient ${GL}_n/O_n$ and a one-parameter family of polar-affine metrics on the manifold ${SPD}_n$. This is stated in the following theorems.



\begin{theorem}[Polar-affine is a square deformation of affine-invariant]
\begin{enumerate}
    \item There exists a one-parameter family of affine-invariant metrics on the quotient ${GL}_n/O_n$.
    \item This family is in bijection with the one-parameter family of affine-invariant metrics on the manifold of SPD matrices thanks to the diffeomorphism $\tau^1:A.O_n\lmto AA^\top$. The corresponding action is $\eta^1:(A,\Sigma)\lmto A\Sigma A^\top$.
    \item This family is also in bijection with a one-parameter family of polar-affine metrics on the manifold of SPD matrices thanks to the diffeomorphism $\tau^2:A.O_n\lmto\sqrt{AA^\top}$. The corresponding action is $\eta^2:(A,\Sigma)\lmto\sqrt{A\Sigma^2A^\top}$.
    \item The diffeomorphism $\pow_2:\fun{({SPD}_n,4g^2)}{({SPD}_n,g^1)}{\Sigma}{\Sigma^2}$ is an isometry between polar-affine metrics $g^2$ and affine-invariant metrics $g^1$.
\end{enumerate}
\end{theorem}

In other words, performing statistical analyses (e.g. a principal component analysis) with the polar-affine metric on covariance matrices is equivalent to performing these statistical analyses with the classical affine-invariant metric on the \textit{square} of our covariance matrix dataset.\\

All the metrics mentioned in Theorem 1 endow their respective space with a structure of a Riemannian symmetric manifold. We recall the definition of that geometric structure and we give the formal statement.

\begin{definition}[Symmetric manifold, Riemannian symmetric manifold]
A manifold $\M$ is symmetric if it is endowed with a family of involutions $(s_x)_{x\in\M}$ called symmetries such that $s_x\circ s_y\circ s_x=s_{s_x(y)}$ and $x$ is an isolated fixed point of $s_x$. It implies that $T_xs_x=-\Id_{T_x\M}$. A Riemannian manifold $(\M,g)$ is symmetric if it is endowed with a family of symmetries that are isometries of $\M$, i.e. that preserve the metric: $g_{s_x(y)}(T_ys_x(v),T_ys_x(w))=g_y(v,w)$ for $v,w\in T_y\M$.
\end{definition}

\begin{theorem}[Riemannian symmetric structure on ${SPD}_n$]
The Riemannian manifold $({SPD}_n,g^1)$, where $g^1$ is an affine-invariant metric, is a Riemannian symmetric space with symmetry $s_\Sigma:\Lambda\lmto\Sigma\Lambda^{-1}\Sigma$. 
The Riemannian manifold $({SPD}_n,g^2)$, where $g^2$ is a polar-affine metric, is also a Riemannian symmetric space whose symmetry is $s_\Sigma:\Lambda\lmto\sqrt{\Sigma^2\Lambda^{-2}\Sigma^2}$.
\end{theorem}


This square deformation of affine-invariant metrics can be generalized into a power deformation to build a family of affine-invariant metrics that we call power-affine metrics. It can even be generalized into any diffeomorphic deformation of SPD matrices. We now develop these families of affine-invariant metrics.




\section{Families of affine-invariant metrics}

There is a theoretical interest in building families comprising some of the known metrics on SPD matrices to understand how one can be deformed into another. For example, power-Euclidean metrics \cite{Dryden10} comprise the Euclidean metric and tends to the log-Euclidean metric \cite{Fillard07} when the power tends to 0. We recall that the log-Euclidean metric is the pullback of the Euclidean metric on symmetric matrices by the symmetric matrix logarithm $\log:{SPD}_n\lto\Sym_n$. There is also a practical interest in defining families of metrics: for example, it is possible to optimize the power to better fit the data with a certain distribution \cite{Dryden10}.

First, we generalize the square deformation by deforming the affine-invariant metrics with a power function $\pow_\theta:\Sigma\in{SPD}_n\lmto\Sigma^\theta=\exp(\theta\log\Sigma)$ to define the power-affine metrics. Then we deform the affine-invariant metrics by any diffeomorphism $f:{SPD}_n\lto{SPD}_n$ to define the deformed-affine metrics.

\subsection{The two-parameter family of power-affine metrics}

We recall that $\M={SPD}_n$ is the manifold of SPD matrices. For a power $\theta\ne 0$, we define the $\theta$-power-affine metric $g^\theta$ as the pullback by the diffeomorphism $\pow_\theta:\Sigma\lmto\Sigma^\theta$ of the affine-invariant metric, scaled by a factor ${1}/{\theta^2}$.

Equivalently, the $\theta$-power-affine metric is the metric invariant under the $\theta$-affine action $\eta^\theta:(A,\Sigma)\lmto(A\Sigma^\theta A^\top)^{1/\theta}$ whose scalar product at $I_n$ coincides with the scalar product $g^1_{I_n}:(V,W)\lmto\alpha\,\tr(VW)+\beta\,\tr(V)\tr(W)$. The $\theta$-affine action induces an isomorphism $V\in T_\Sigma M\lmto V_\theta=\frac{1}{\theta}\Sigma^{-\theta/2}\,\partial_V\pow_\theta(\Sigma)\,\Sigma^{-\theta/2}\in T_{I_n}\M$ between tangent spaces. The $\theta$-power-affine metric is given by:
\begin{equation}
    g^\theta_\Sigma(V,W)=\alpha\,\tr(V_\theta W_\theta)+\beta\,\tr(V_\theta)\tr(W_\theta).
\end{equation}

Because a scaling factor is of low importance, we can set $\alpha=1$ and consider that this family is a two-parameter family indexed by $\beta>-{1}/{n}$ and $\theta\ne 0$.

We have chosen to define the metric $g^\theta$ so that the power function $\pow_\theta:(\M,\theta^2g^\theta)\lto(\M,g^1)$ is an isometry. Why this factor $\theta^2$? The first reason is for consistence with previous works: the analogous power-Euclidean metrics have been defined with that scaling \cite{Dryden10}. The second reason is for continuity: when the power tends to 0, the power-affine metric tends to the log-Euclidean metric.

\begin{theorem}[Power-affine tends to log-Euclidean for $\theta\rightarrow 0$]
Let $\Sigma\in\M$ and $V,W\in T_\Sigma\M$. 
Then $\lim_{\theta\rightarrow 0}{g^\theta_\Sigma(V,W)}=g^{LE}_\Sigma(V,W)$
where the log-Euclidean metric is $g^{LE}_\Sigma(V,W)=\alpha\, \tr(\partial_V\!\log(\Sigma)\,\partial_W\!\log(\Sigma))+\beta\,\tr(\partial_V\!\log(\Sigma))\tr(\partial_W\!\log(\Sigma))$.
\end{theorem}

\subsection{The continuum of deformed-affine metrics}
In the following, we call  a diffeomorphism $f:{SPD}_n\lto{SPD}_n$ a deformation. We define the $f$-deformed-affine metric $g^f$ as the pullback by the diffeomorphism $f$ of the affine-invariant metric, so that $f:(\M,g^f)\lto(\M,g^1)$ is an isometry. (Regarding the discussion before the Theorem 3, $g^{\pow_\theta}=\theta^2g^\theta$.)

The $f$-deformed-affine metric is invariant under the $f$-affine action $\eta^f:(A,\Sigma)\lmto f^{-1}(Af(\Sigma)A^\top)$. It is given by $g^f_\Sigma(V,W)=\alpha\tr(V_fW_f)+\beta\tr(V_f)\tr(W_f)$
where $V_f=f(\Sigma)^{-1/2}\, \partial_Vf(\Sigma) \, f(\Sigma)^{-1/2}$. The basic Riemannian operations are obtained by pulling back the affine-invariant operations.

\begin{theorem}[Basic Riemannian operations]
For SPD matrices $\Sigma,\Lambda\in\M$ and a tangent vector $V\in T_\Sigma\M$, we have at all time $t\in\R$:
\begin{center}
$\begin{array}{|c|l|}
\hline
\mathrm{Geodesics} & \gamma^f_{(\Sigma,V)}(t)=f^{-1}(f(\Sigma)^{1/2}\exp(tf(\Sigma)^{-1/2}T_\Sigma f(V)f(\Sigma)^{-1/2})f(\Sigma)^{1/2}) \\
\hline
\mathrm{Logarithm} & \mathrm{Log}^f_\Sigma(\Lambda)=(T_\Sigma f)^{-1}(f(\Sigma)^{1/2}\log(f(\Sigma)^{-1/2}f(\Lambda)f(\Sigma)^{-1/2})f(\Sigma)^{1/2})\\
\hline
\mathrm{Distance} & d_f(\Sigma,\Lambda)=d_1(f(\Sigma),f(\Lambda))=\sum_{k=1}^n{(\log\lambda_k)^2}\\
\hline
\end{array}$
\end{center}
where $\lambda_1,...,\lambda_n$ are the eigenvalues of the symmetric matrix $f(\Sigma)^{-1/2}f(\Lambda)f(\Sigma)^{-1/2}$.
\end{theorem}

All tensors are modified thanks to the pushforward $f_*$ and pullback $f^*$ operators, e.g. the Riemann tensor of the $f$-deformed metric is $R^f(X,Y)Z=f^*(R(f_*X,f_*Y)(f_*Z))$. As a consequence, the deformation $f$ does not affect the values taken by the sectional curvature and these metrics are negatively curved.


From a computational point of view, it is very interesting to notice that  the identification $\mc{L}'_\Sigma:V\in T_\Sigma\M\lmto V'=T_\Sigma f(V)\in T_{f(\Sigma)}\M$ simplifies the above expressions by removing the differential $T_\Sigma f$. This change of basis can prevent from numerical approximations of the differential but one must keep in mind that $V\ne V'$ in general. This identification was already used for the polar-affine metric ($f=\pow_2$) in \cite{Zhang18} without explicitly mentioning.

\section{Interesting subfamilies of deformed-affine metrics}

Some deformations have already been used in applications. For example, the family $A_r:\diag(\lambda_1,\lambda_2,\lambda_3)\lmto\diag(a_1(r)\lambda_1,a_2(r)\lambda_2,a_3(r)\lambda_3)$ where $\lambda_1\gs\lambda_2\gs\lambda_3>0$ was proposed to map the anisotropy of water measured by diffusion tensors to the one of the diffusion of tumor cells in tumor growth modeling \cite{Jbabdi05}. The inverse function $\inv=\pow_{-1}:\Sigma\lmto\Sigma^{-1}$ or the adjugate function $\adj:\Sigma\lmto\det(\Sigma)\Sigma^{-1}$ were also proposed in the context of DTI \cite{Lenglet04,Fuster16}. Let us find some properties satisfied by some of these examples. We define the following subsets of the set $\mc{F}=\Diff({SPD}_n)$ of diffeomorphisms of ${SPD}_n$.

(Spectral) $\mc{S}=\{f\in\mc{F}|\forall U\in O_n,\forall D\in\Diag^{++}_n,f(UDU^\top)=Uf(D)U^\top\}$.
Spectral deformations are characterized by their values on sorted diagonal matrices so the deformations described above are spectral: $A_r,\adj,\pow_\theta\in\mc{S}$.\\
For a spectral deformation $f\in\mc{S}$, $f(\R_+^*I_n)=\R_+^*I_n$ so we can unically define a smooth diffeomorphism $f_0:\R_+^*\lto\R_+^*$ by $f(\lambda I_n)=f_0(\lambda)I_n$.
    
(Univariate) $\mc{U}=\{f\in\mc{S}|f(\diag(\lambda_1,...,\lambda_n))=\diag(f_0(\lambda_1),...,f_0(\lambda_n))\}$.
The power functions are univariate. Any polynomial $P=\lambda X\prod_{k=1}^p(X-a_i)$ null at 0, with non-positive roots $a_i\ls 0$ and positive coefficient $\lambda>0$, also gives rise to a univariate deformation.
    
(Diagonally-stable) $\mc{D}=\{f\in\mc{F}|f(\Diag^{++}_n)\subset\Diag^{++}_n\}$.
The deformations described above $A_r,\adj,\pow_\theta$ and the univariate deformations are clearly diagonally-stable: $A_r,\adj,\pow_\theta\in\mc{D}$ and $\mc{U}\subset\mc{D}\cap\mc{S}$.

(Log-linear) $\mc{L}=\{f\in\mc{F}|\log_*f=\log\circ\, f\circ\exp\,\,\text{is linear}\}$.
The adjugate function and the power functions are log-linear deformations. More generally, the functions $f_{\lambda,\mu}:\Sigma\lmto(\det\Sigma)^{\frac{\lambda-\mu}{n}}\Sigma^\mu$ for $\lambda,\mu\ne0$, are log-linear deformations. We can notice that the $f_{\lambda,\mu}$-deformed-affine metric belongs to the one-parameter family of $\mu$-power-affine metrics with $\beta=\frac{\lambda^2-\mu^2}{n\mu^2}>-\frac{1}{n}$.

The deformations $f_{\lambda,\mu}$ just introduced are also spectral and the following result states that they are the only spectral log-linear deformations.

\begin{theorem}[Characterization of the power-affine metrics]
If $f\in\mc{S}\cap\mc{L}$ is a spectral log-linear diffeomorphism, then there exist real numbers $\lambda,\mu\in\R^*$ such that $f=f_{\lambda,\mu}$ and the $f$-deformed-affine metric is a $\mu$-power-affine metric.
\end{theorem}

The interest of this theorem comes from the fact that the group of spectral deformations and the vector space of log-linear deformations have large dimensions while their intersection is reduced to a two-parameter family. This strong result is a consequence of the theory of Lie group representations because the combination of the spectral property and the linearity makes $\log_*f$ a homomorphism of $O_n$-modules (see the sketch of proof below).\\

\textit{Sketch of the proof.}
Thanks to Lie group representation theory, the linear map $F=\log_*f:{Sym}_n\lto{Sym}_n$ appears as a homomorphism of $O_n$-modules for the representation $\rho:P\in O_n\lmto(V\lmto PVP^\top)\in GL({Sym}_n)$. Once shown that ${Sym}_n=\Vect(I_n)\oplus\ker{\tr}$ is a $\rho$-irreducible decomposition of ${Sym}_n$ and that each one is stable by $F$, then according to Schur's lemma, $F$ is homothetic on each subspace, i.e. there exist $\lambda,\mu\in\R^*$ such that for $V\in\Sym_n$, $F(V)=\lambda\frac{\tr(V)}{n}I_n+\mu\left(V-\frac{\tr(V)}{n}I_n\right)=\log_*f_{\lambda,\mu}(V)$, so $f=f_{\lambda,\mu}$.\\

\section{Conclusion}
We have shown that the polar-affine metric is a square deformation of the affine-invariant metric and this process can be generalized to any power function or any diffeomorphism on SPD matrices. It results that the invariance principle of symmetry is not sufficient to distinguish all these metrics, so we should find other principles to limit the scope of acceptable metrics in statistical computing. We have proposed a few characteristics (spectral, diagonally-stable, univariate, log-linear) that include some functions on tensors previously introduced. 
Future work will focus on studying the effect of such deformations on real data and on extending this family of metrics to positive semi-definite matrices. Finding families that comprise two non-cousin metrics could also help understand the differences between them and bring principles to make choices in applications.

\paragraph{Acknowledgements.} This project has received funding from the European Research Council (ERC) under the European Union’s Horizon 2020 research and innovation program (grant G-Statistics agreement No 786854). This work has been supported by the French government, through the UCAJEDI Investments in the Future project managed by the National Research Agency (ANR) with the reference number ANR-15-IDEX-01.

\bibliographystyle{unsrt}
\bibliography{biblio}

\end{document}